\documentclass[12pt,a4]{amsart}
\usepackage[a4paper, left=28mm, right=28mm, top=28mm, bottom=34mm]{geometry}
\usepackage[all]{xy}
\usepackage{amsmath}
\usepackage{amssymb}
\usepackage{amsthm}
\usepackage{listings}
\usepackage{mathrsfs}
\theoremstyle{definition}
\newtheorem{thm}{Theorem}[section]
\newtheorem{lem}[thm]{Lemma}

\newtheorem{prop}[thm]{Proposition}

\theoremstyle{definition}

\newtheorem{rmk}[thm]{Remark}

\newtheorem{cond}[thm]{Condition}
\newtheorem{exa}[thm]{Example}
\numberwithin{equation}{section}
\def\F{{\mathbb F}}
\def\P{{\mathbb P}}
\def\Q{{\mathbb Q}}

\def\Arf{\mathop{\mathrm{Arf}}\nolimits}
\def\Gal{\mathop{\mathrm{Gal}}\nolimits}
\def\Jac{\mathop{\mathrm{Jac}}\nolimits}
\def\Sp{\mathop{\mathrm{Sp}}\nolimits}
\def\chara{\mathop{\mathrm{char}}}
\def\pr{\text{\rm pr}}
\def\sep{\text{\rm sep}}
\begin{document}

\title[Local-global property for bitangents]
{The local-global property for bitangents of plane quartics}

\author[Y.\ Ishitsuka]{Yasuhiro Ishitsuka}
\address{Center for Science Adventure and Collaborative Research Advancement, Graduate School of Science, Kyoto University, Kyoto 606-8502, Japan}
\address{Department of Mathematics, Faculty of Science, Kyoto University, Kyoto 606-8502, Japan}
\email{yasu-ishi@math.kyoto-u.ac.jp}

\author[T.\ Ito]{Tetsushi Ito}
\address{Department of Mathematics, Faculty of Science, Kyoto University, Kyoto 606-8502, Japan}
\address{Mathematical Science Team, RIKEN Center for Advanced Intelligence Project (AIP), 1-4-1 Nihonbashi Chuo-ku, Tokyo 103-0027, Japan}
\email{tetsushi@math.kyoto-u.ac.jp}

\author[T.\ Ohshita]{Tatsuya Ohshita}
\address{Department of Mathematics, Faculty of Science and Technology, Keio University, 3-14-1 Hiyoshi, Kohoku-ku Yokohama-shi, Kanagawa 223-8522, Japan}
\email{ohshita@math.keio.ac.jp}

\author[T.\ Taniguchi]{Takashi Taniguchi}
\address{Department of Mathematics, Graduate School of Science, Kobe University, Kobe 657-8501, Japan}
\email{tani@math.kobe-u.ac.jp}

\author[Y.\ Uchida]{Yukihiro Uchida}
\address{Department of Mathematical Sciences, Graduate School of Science, Tokyo Metropolitan University, 1-1 Minami-Osawa, Hachioji, Tokyo 192-0397, Japan}
\email{yuchida@tmu.ac.jp}

\date{April 20, 2020}


\maketitle

\begin{abstract}
We study the arithmetic of bitangents of smooth quartics over global fields.
With the aid of computer algebra systems and using Elsenhans--Jahnel's results on the inverse Galois problem for bitangents,
we show that, over any global field of characteristic different from $2$,
there exist smooth quartics
which have bitangents over every local field,
but do not have bitangents over the global field.
We give an algorithm to find such quartics explicitly,
and give an example over $\Q$.
We also discuss a similar problem concerning symmetric determinantal representations.
This paper is a summary of the first author's talk at the JSIAM JANT workshop
on algorithmic number theory in March 2019. Details will appear elsewhere.
\end{abstract}

\section{Introduction}

Let $C \subset \P^2$ be a smooth quartic
over a field $k$ of characteristic different from $2$.
It is defined by a homogeneous polynomial
$f(X,Y,Z)$ of degree $4$:
\[
C = \{ [X : Y : Z] \in \P^2 \mid f(X,Y,Z) = 0 \}.
\]
By B\'ezout's theorem,
the intersection of $C$ with a line $L \subset \P^2$ consists of four points,
counted with their multiplicities.
A line $L \subset \P^2$ is called
a \textit{bitangent} of $C$
if one of the following conditions is satisfied:
\begin{itemize}
    \item $L$ tangents to $C$ at two distinct points, or
    \item $L$ tangents quadruply to $C$ at one point.
    (In this case, the line $L$ is also called a \textit{hyperflex line} of $C$.)
\end{itemize}
It is known that every smooth quartic has exactly $28$ bitangents
over an algebraic closure of $k$ \cite[Chapter 6]{Dolgachev}.

Bitangents play an important role in the study of the arithmetic of quartics.
For example, Bruin--Poonen--Stoll gave an application to the calculation of the Mordell--Weil group and the rational points \cite{BruinPoonenStoll}.
However, since it involves computations on number fields of degree 28 in general, 
it remains a difficulty to treat the arithmetic of bitangents explicitly.

In this paper, we shall consider the problem whether a given smooth quartic over a global field $K$ has a bitangent over $K$.
Concerning this problem, a natural theoretical question is whether the bitangents satisfy the \textit{local-global property} (or \textit{Hasse principle}) or not.
The aim of this paper is to answer this question \textit{negatively}
with the aid of computer algebra systems.

Here is the main theorem of this paper.

\begin{thm}
\label{MainTheorem}
Let $K$ be a global field of characteristic different from $2$.
(It is a finite extension of $\Q$ or $\F_p(T)$ with $p \neq 2$.)
Then, there exists a smooth quartic $C \subset \P^2$ over $K$ such that
\begin{itemize}
\item $C$ has a bitangent over $K_v$ for every place $v$ of $K$ (including the infinite places when $K$ is a number field), but
\item $C$ does not have a bitangent over $K$.
\end{itemize}
\end{thm}

Our proof of Theorem \ref{MainTheorem}
is a combination of certain group theoretic results on $\Sp_6(\F_2)$
obtained with the aid of computer algebra systems
and the results on the inverse Galois problem for bitangents
recently obtained by Elsenhans--Jahnel
\cite{ElsenhansJahnel:CayleyOctad},
\cite{ElsenhansJahnel:SteinerHexad}.

\begin{rmk}
After this work was completed, 
Jahnel and Loughran told the authors that it is also possible to
construct smooth quartics satisfying the conditions in
Theorem \ref{MainTheorem} by the results in \cite{JahnelLoughran}.
(See Remark \ref{Remark:JahnelLoughranl}.)
\end{rmk}

\begin{rmk}
In this paper, we shall also give an algorithm to find
quartics failing the local-global property for bitangents explicitly.
(See Section \ref{Sect:Algorithm}.)
For a given global field $K$,
the hardest step is to find a certain Galois extension $L/K$ of degree 32.
It seems difficult to give a sharp estimate of the time complexity because this step is related to the effective Chebotarev density theorem.
For more details, see Appendix B.
\end{rmk}

\section{Bitangents and quadratic forms}

Let $k$ be a field of characteristic different from $2$,
and $k^{\sep}$ a separable closure of $k$.
Let $C \subset \P^2$ be a smooth quartic over $k$.
The Jacobian variety $\Jac(C)$ is an abelian variety of dimension $3$. The group $\Jac(C)[2]$ of $k^{\sep}$-rational points on $\Jac(C)$ killed by $2$ is a $6$-dimensional vector space over $\F_2$.
It is equipped with the action of $\Gal(k^{\sep}/k)$
and the Weil pairing
\[
\langle,\rangle \colon \Jac(C)[2] \times \Jac(C)[2]
\to \{ \pm 1 \} \cong \F_2.
\]

We fix a symplectic basis $\{ e_1,e_2,e_3,f_1,f_2,f_3 \}$
of $\Jac(C)[2]$,
i.e.\ $\langle e_i, e_j \rangle = \langle f_i, f_j \rangle = 0$ and $\langle e_i, f_j \rangle = \delta_{i,j}$
for all $1 \leq i,j \leq 3$.
We have an isomorphism $\Jac(C)[2] \cong \F_2^{\oplus 6}$.
The action of $\Gal(k^{\sep}/k)$ gives a continuous map
\[
\rho_C \colon \Gal(k^{\sep}/k) \to \Sp(\Jac(C)[2], \langle,\rangle) \cong \Sp_6(\F_2).
\]

A map $Q \colon \F_2^{\oplus 6} \to \F_2$ is called
a \textit{quadratic form with polar form $\langle,\rangle$}
if
\[ Q(x+y) - Q(x) - Q(y) = \langle x,y \rangle \]
is satisfied for every $x,y \in \F_2^{\oplus 6}$.
The \textit{Arf invariant} of $Q$ is defined by
\[
\Arf(Q) := \sum_{1 \leq i \leq 3} Q(e_i) Q(f_i) \in \F_2,
\]
which is independent of the choice of a symplectic basis.

Let $\Omega^{+}$ (resp.\ $\Omega^{-}$) be the subset of $\Omega$
consisting of quadratic forms of Arf invariant $0$ (resp.\ $1$).
The set $\Omega^{+}$ (resp.\ $\Omega^{-}$) 
has $36$ (resp.\ $28$) elements.
The symplectic group $\Sp_6(\F_2)$ acts transitively on
$\Omega^{+}$ and $\Omega^{-}$.

\begin{lem}[Mumford]
\label{Lemma:Bijection}
There is a bijection between the set $\Omega^{-}$ and
the set of bitangents of $C$.
\end{lem}

See \cite{Mumford}, \cite[Proposition 6.2]{IshitsukaIto} for details.
Since $\chara(k) \neq 2$, every bitangent of $C$ is defined over $k^{\sep}$.
The bijection in Lemma \ref{Lemma:Bijection}
is equivariant under the action of $\Gal(k^{\sep}/k)$,
where $\Gal(k^{\sep}/k)$ acts on $\Omega^{-}$
through the map
$\rho_C \colon \Gal(k^{\sep}/k) \to \Sp_6(\F_2)$.

\begin{rmk}
For a `generic' choice of $C$, Harris and Shioda proved 
that the map $\rho_C$ is surjective,
at least when $\chara(k) \notin \{ 3, 5, 7, 11, 29, 1229 \}$; see \cite[p.721]{Harris}, \cite[Theorem 7]{Shioda}.
For explicit examples of quartics with surjective $\rho_C$,
see \cite[p.69, Example]{Shioda}, \cite[p.26, Corollary 3]{Erne}.
Shioda also constructed smooth quartics over $\Q$
such that all the $28$ bitangents are defined over $\Q$;
see \cite[(6.6)]{Shioda}.
\end{rmk}

\section{Group theoretic results}

Let $U_{36} \subset \Sp_6(\F_2)$
be the stabilizer of an element of $\Omega^{+}$.
The action of $\Sp_6(\F_2)$ on
$\F^{\oplus 6}_2 \setminus \{0\}$ is transitive.
The stabilizer of a non-zero vector is denoted by
$U_{63} \subset \Sp_6(\F_2)$.
It is known that 
$U_{36}, U_{63}$ are maximal subgroups of $\Sp_6(\F_2)$,
and every subgroup of $\Sp_6(\F_2)$
of index $36$, $63$ is conjugate to
$U_{36}$, $U_{63}$, respectively.

We shall consider the following condition.

\begin{cond}
\label{Condition:LocalGlobal}
Let $G \subset \Sp_6(\F_2)$ be a subgroup.
We say $G$ satisfies the {\em condition} $(\ast)^{+}$ (resp.\ $(\ast)^{-}$)
if the following conditions are satisfied:
\begin{itemize}
\item No element of $\Omega^{+}$ (resp.\ $\Omega^{-})$ is
fixed by every element of $G$.
\item For every $g \in G$, the action of $g$ on
$\Omega^{+}$ (resp.\ $\Omega^{-}$)
has a fixed point.
\end{itemize}
\end{cond}

The following results can be confirmed by GAP.
(For a sample source code for GAP, see Appendix A.)

\begin{prop}
\label{Proposition:FixedPointProperties}
\begin{enumerate}
\item $U_{36}$ has $296$ subgroups, up to conjugation.
Among them, $35$ subgroups satisfy the condition $(\ast)^{-}$;
all the $35$ subgroups are solvable.

\item $U_{63}$ has $1916$ subgroups, up to conjugation.
Among them, $548$ subgroups satisfy both of the condition $(\ast)^{+}$ and $(\ast)^{-}$;
$536$ of the $548$ subgroups are solvable.
\end{enumerate}
\end{prop}

\begin{prop}
\label{Proposition:FixedPointProperties2}
$\Sp_6(\F_2)$ has $6$ subgroups isomorphic to
$\F_2^{\oplus 5}$, up to conjugation.
All of them satisfy both of the conditions
$(\ast)^{+}$ and $(\ast)^{-}$.
\end{prop}

\begin{rmk}
Any subgroup of $U_{36}$ does not satisfy the condition $(\ast)^{+}$
since, by definition, it fixes at least one element of
$\Omega^{+}$.
\end{rmk}


\section{The inverse Galois problems}

The inverse Galois problem asks
whether every finite group is realized
as the Galois group of a number field over $\Q$.
It is an open problem in general.
It was proved by Shafarevich for solvable groups.
Sonn observed that Shafarevich's solution to
the inverse Galois problem yields
a Galois extension such that
every decomposition group is cyclic;
see \cite[Theorem 2]{Sonn} for details.

\begin{thm}[Shafarevich, Sonn]
\label{Theorem:ShafarevichSonn}
Let $G$ be a finite solvable group, and $K$ a global field.
Then there exists a finite Galois extension $L/K$
such that
$\Gal(L/K)$ is isomorphic to $G$,
and that, for every  place $v$ of $K$,
the decomposition group of $L/K$ at $v$ is cyclic.
\end{thm}

Recently, Elsenhans--Jahnel studied
an analogue of the inverse Galois problem
for bitangents of quartics.
Let $k$ be a field of characteristic different from $2$, and
\[ \rho \colon \Gal(k^{\sep}/k) \to \Sp_6(\F_2) \]
a continuous homomorphism.
They asked whether $\rho$
is realized as the map $\rho_C$ associated with a smooth quartic $C \subset \P^2$ over $k$,
up to conjugation by an element of $\Sp_6(\F_2)$.
The following result is proved in
\cite{ElsenhansJahnel:CayleyOctad}
(resp.\ \cite{ElsenhansJahnel:SteinerHexad})
when the image of $\rho$
is contained in a conjugate of $U_{36}$
(resp.\ $U_{63}$).

\begin{thm}[Elsenhans--Jahnel]
\label{ElsenhansJahnel}
If the image of $\rho$
is contained in a conjugate of $U_{36}$ or $U_{63}$,
there exists a smooth quartic $C \subset \P^2$
over $k$ such that the maps $\rho, \rho_C$
are conjugate by an element of $\Sp_6(\F_2)$.
\end{thm}


\section{Proof of Theorem \ref{MainTheorem}}

Theorem \ref{MainTheorem} is proved by combining 
Lemma \ref{Lemma:Bijection},
Proposition \ref{Proposition:FixedPointProperties},
Theorem \ref{Theorem:ShafarevichSonn},
and Theorem \ref{ElsenhansJahnel}.
\begin{enumerate}
\item Take a solvable subgroup
\[ G \subset \Sp_6(\F_2) \]
contained in $U_{36}$ or $U_{63}$ which satisfies the condition $(\ast)^{-}$.
Such a subgroup exists by Proposition \ref{Proposition:FixedPointProperties}.

\item Take a finite Galois extension $L/K$ such that
\[ \Gal(L/K) \cong G \]
and, for every place $v$ of $K$,
the decomposition group of $L/K$ at $v$ is cyclic.
Such an extension exists by Theorem \ref{Theorem:ShafarevichSonn}.

\item Consider the composite of the following maps:
\[
\rho \colon \Gal(K^{\sep}/K) \to \Gal(L/K) \cong G \hookrightarrow \Sp_6(\F_2).
\]
Take a smooth quartic $C \subset \P^2$ over $K$
such that the maps $\rho, \rho_C$
are conjugate by an element of $\Sp_6(\F_2)$.
Such a quartic exists by Theorem \ref{ElsenhansJahnel}.

\item By Lemma \ref{Lemma:Bijection}, there is
a $\Gal(K^{\sep}/K)$-equivariant bijection
between the set $\Omega^{-}$ and
the set of bitangents of $C$ over $K^{\sep}$.
For a place $v$ of $K$,
$\Gal(K_v^{\sep}/K_v)$ is embedded into $\Gal(K^{\sep}/K)$.
The image
\[ G_v := \rho(\Gal(K_v^{\sep}/K_v)) \]
is a cyclic group by (2).
Let $g_v$ be a generator of $G_v$.
Since $G$ satisfies the condition $(\ast)^{-}$ by (1),
there is an element of $\Omega^{-}$ fixed by $g_v$.
Hence the quartic $C$ has a bitangent over $K_v$.
Moreover, since no element of $\Omega^{-}$ is fixed by every element of $G$, the quartic $C$ does not have a bitangent over $K$.
Therefore, the quartic $C$ satisfies all the conditions of
Theorem \ref{MainTheorem}.
\end{enumerate}

\section{Description of algorithm}
\label{Sect:Algorithm}

Here we shall give an algorithm obtaining
smooth quartics failing the local-global property
for bitangents.

We take a subgroup
$G \subset U_{63} \subset \Sp_6(\F_2)$
isomorphic to $\F_2^{\oplus 5}$.
We can use Proposition \ref{Proposition:FixedPointProperties2}
for such a subgroup.
Explicitly, the group $U_{63}$ is isomorphic to
$(\F_2 \wr S_6) \cap A_{12}$;
see \cite[Corollary 2.19]{ElsenhansJahnel:SteinerHexad}.
The wreath product $\F_2 \wr S_6$ contains
a subgroup isomorphic to $\F_2^{\oplus 6}$.
Let $G$ be the kernel of the sum $\F_2^{\oplus 6} \to \F_2$.
Then $G$ is isomorphic to $\F_2^{\oplus 5}$,
and we have an embedding
\[ G \cong \F_2^{\oplus 5} \hookrightarrow 
(\F_2 \wr S_6) \cap A_{12} \cong U_{63}. \]

Our algorithm consists of $3$ steps.

\subsection*{Step 1 (Find a Galois extension $L/K$.)}

Take a Galois extension $L/K$ such that
\[ \Gal(L/K) \cong \F_2^{\oplus 5} \]
and that every decomposition group is cyclic.
(The existence of such $L/K$ is guaranteed by Theorem \ref{Theorem:ShafarevichSonn}.)

Assume that $L$ is described as 
\[
    L = K \left( \sqrt{a_1},\,\sqrt{a_2},\,\sqrt{a_3},\,\sqrt{a_4},\,\sqrt{a_5}\right).
\]
for some $a_1,a_2,a_3,a_4,a_5 \in K$.
We put $a_6 = a_1 a_2 a_3 a_4 a_5 u^2$ for some $u \in K^{\times}$ and set
\[
    F(S, T) = (a_1 S-T)(a_2S-T) \cdots (a_6S-T).
\]
By construction, the splitting field of
$F(1,T^2)$ is $L$ and $F(1,0) = (a_1 a_2 \cdots a_5u)^2$.

\subsection*{Step 2 (Construct a conic bundle $B \subset \P^1 \times \P^2$.)}

There exists a unique pair of binary quadratic forms $g(S,T), h(S,T)$ satisfying
\[ \det M(S,T) = -F(S, T) \quad \mathrm{and} \quad g(1,0) = a_1 a_2 a_3 a_4 a_5 u, \]
where
\[
    M(S,T) = \begin{pmatrix}
        -ST + T^2 & ST  & g(S,T) \\
        ST        & S^2 & T^2    \\
        g(S,T)    & T^2 & h(S,T)
    \end{pmatrix}.
\]

Let $B \subset \P^1 \times \P^2$ be the hypersurface defined by
\[
    \begin{pmatrix}
        X & Y & Z
    \end{pmatrix}
    M(S,T)
    \begin{pmatrix}
        X \\ Y \\ Z
    \end{pmatrix}
    = 0.
\]
Here $[S : T]$ is the projective coordinate of $\P^1$ and 
$[X : Y : Z]$ is that of $\P^2$.
The first projection $\pr_1 \colon B \to \P^1$ defines
a conic bundle structure with
six degenerate fibers with prescribed Galois action;
see \cite[Proposition 3.5]{ElsenhansJahnel:SteinerHexad}.

\subsection*{Step 3 (Calculate the branching locus.)}

The composite of the embedding
\[ \iota \colon B \hookrightarrow \P^1 \times \P^2 \]
and the projection
\[ \pr_2 \colon \P^1 \times \P^2 \to \P^2 \]
is a double cover.
Its ramification locus is the desired quartic $C \subset \P^2$.
(If $C$ is not smooth, then take other parameters
$a_1, a_2, \ldots, a_5$ and $u$,
and calculate the quartic $C$ again.
For a `generic' choice of $a_1, a_2, \ldots, a_5$ and $u$,
the quartic $C$ is smooth; see \cite[Proposition 3.5]{ElsenhansJahnel:SteinerHexad}.)

Applying the above algorithm, we can construct a smooth quartic
$C \subset \P^2$
over $K$ such that the map $\rho_C$ is
conjugate to the composite of the following maps:
\[
\Gal(K^{\sep}/K) \to \Gal(L/K) \cong G \hookrightarrow \Sp_6(\F_2).
\]
The smooth quartic $C \subset \P^2$ satisfies
all the conditions of Theorem \ref{MainTheorem}
by Proposition \ref{Proposition:FixedPointProperties2}.

\begin{rmk}[Jahnel--Loughran]
\label{Remark:JahnelLoughranl}
Here is another construction of smooth quartics
satisfying the conditions of Theorem \ref{MainTheorem}
by the results in \cite{JahnelLoughran}.
Take two closed points $P,Q$ of the projective plane $\P^2$ over $K$ (as a $K$-scheme) of degree $3,4$, respectively, such that
\begin{itemize}
\item the union $P \sqcup Q$ lies in general position, and
\item the blow-up of $\P^2$ along $Q$ is a del Pezzo surface of degree $5$ which does not satisfy the local-global property for lines; see \cite[Section 3.6]{JahnelLoughran}.
\end{itemize}
Let $X$ be the blow-up of $\P^2$ along $P \sqcup Q$.
The branching locus $C \subset \P^2$ of the anticanonical map $\pi \colon X \to \P^2$
is a smooth quartic.
We have a $2:1$ map from the $56$ lines on $X$ to the $28$ bitangents of $C$.
For every place $v$ of $K$, $X$ has a line over $K_v$
by \cite[Lemma 3.9]{JahnelLoughran}.
Hence $C$ has a bitangent over $K_v$.
Moreover, $C$ does not have a bitangent over $K$.
In fact, if $L \subset \P^2$ were a bitangent over $K$,
its inverse image $\pi^{-1}(L) = L_1 \cup L_2$ is a union of two lines on $X$.
The lines on $X$ are classified into the types (i)-(iv)
in \cite[Remark 2.8]{ElsenhansJahnel:SteinerHexad}.
Since $L_1, L_2$ have different types
and the type of a line is Galois invariant,
both of $L_1, L_2$ are defined over $K$.
It contradicts \cite[Lemma 3.8]{JahnelLoughran}.
\end{rmk}

\section{An example}
\label{Section:Example}

We consider the case $K = \Q$.
We put
\[ b_1 = -1,\quad b_2 = 17,\quad b_3 = 89,\quad b_4 = 257,\quad b_5 = 769, \]
and put $K_i = \Q(\sqrt{b_i})$ for every $1 \leq i \leq 5$.
Then, for every $1 \leq i \leq 5$, only one prime number is ramified in $K_i/\Q$.
For every $1 \leq i,j \leq 5$ with $i \neq j$,
the prime number $p_i$ splits in $K_j/\Q$.
(Here $p_i$ is the unique prime number ramified in $K_i/\Q$.)
We put
\[ L = K_1 K_2 K_3 K_4 K_5, \]
and
\[ a_1 = b_1b_5, \quad a_2 = b_2b_4, \quad a_3 = b_3, \quad a_4 = b_4, \quad a_5 = b_5, \quad u = -b_4^{-1}b_5^{-1}. \]

Then $L/K$ satisfies the conditions in Step 1.
We have
\begin{align*}
F(S, T) &= (-769S - T)(4369S - T)(89S - T) \\
&\quad (257S - T)(769 S- T)(-1513 S - T), \\
g(S,T) &=  \frac{1}{8} (2392149832S^2 + 35008837ST + 12804T^2),\\ 
h(S,T) &= -\frac{1}{64}(251582881045706064S^2 \\
&\quad + 1084638148302617ST + 594847875240T^2).
\end{align*}

The quartic is computed as
\begin{align*}
& \quad 4096 X^{4} - 16384 X^{3} Y - 9869943810048 X^{3} Z \\ 
& + 143396196352 X^{2} Y Z - 52445184 X Y^{2} Z - 32768 Y^{3} Z \\
& + 64826445425191482752 X^{2} Z^{2} \\
& - 277686962456893696 X Y Z^{2} \\
& + 152281056061440 Y^{2} Z^{2} \\
& - 917870567374331469445024 X Z^{3} \\
& + 128810435095401504768 Y Z^{3} \\
& + 577825743806146102974275227249 Z^{4} = 0.
\end{align*}

It gives an example of smooth quartics over $\Q$
failing the local-global property for bitangents.


\section{The local-global property for symmetric determinantal representations}

We say a smooth quartic $C \subset \P^2$ over a field $k$
admits a \textit{symmetric determinantal representation} over $k$
if there exist symmetric matrices
$M_1,M_2,M_3 \in \mathrm{Mat}_{4}(k)$
of size $4 \times 4$ such that the equation
\[ \det(X M_1 + Y M_2 + Z M_3) = 0 \]
defines the quartic $C \subset \P^2$.
(For explicit examples, see \cite{IIO}.)

By the same method as above, it is possible to obtain
smooth quartics failing the local-global property for
symmetric determinantal representations.

\begin{thm}
\label{MainTheorem:SDR}
For any global field $K$ of characteristic different from $2$,
there exists a smooth quartic $C \subset \P^2$ over $K$
such that
\begin{itemize}
\item $C$ admits a symmetric determinantal representation over $K_v$ for every place $v$ of $K$, but
\item $C$ does not admit a symmetric determinantal representation over $K$.
\end{itemize}
\end{thm}

Here is a sketch of the proof.
For a smooth quartic $C \subset \P^2$
with a $K$-rational point,
it admits a symmetric determinantal representation over $K$
if and only if 
there exists a $\Gal(K^{\sep}/K)$-invariant quadratic form in $\Omega^{+}$; see \cite[Theorem 2.2, Corollary 6.3]{IshitsukaIto}.
Taking a subgroup $G \subset U_{63}$
satisfying the condition $(\ast)^{+}$,
we find smooth quartics over $K$
satisfying the conditions of Theorem \ref{MainTheorem:SDR}
by the same way as in the case of bitangents.

\begin{exa}
Quartics constructed by the algorithm in Section \ref{Sect:Algorithm}
always have the $\Q$-rational point $[0:1:0]$.
By Proposition \ref{Proposition:FixedPointProperties2},
the example of quartic in Section \ref{Section:Example}
satisfies the conditions of Theorem \ref{MainTheorem:SDR}.
(But the quartics constructed by the method described in
Remark \ref{Remark:JahnelLoughranl} do not satisfy these conditions.)
\end{exa}

\begin{rmk}
In \cite{IshitsukaIto},
it is proved that smooth quartics over number fields do not satisfy the local-global property for symmetric determinantal representations.
The quartics constructed in \cite{IshitsukaIto}
are defined over number fields of large degree.
For quartics over $\Q$,
this problem was stated in \cite[Problem 1.6 (1)]{IshitsukaIto},
but not answered there.
\end{rmk}

\begin{rmk}
Theorem \ref{MainTheorem:SDR} does not hold in characteristic $2$.
In fact, over a global field of characteristic $2$,
any smooth plane curve of any degree satisfies
the local-global property for symmetric determinantal representations;
see \cite{IshitsukaIto:Char2} for details.
\end{rmk}

\appendix
\section{Sample source codes for GAP and SageMath}
\label{Appendix:SourceCode}

We proved
Proposition \ref{Proposition:FixedPointProperties} and
Proposition \ref{Proposition:FixedPointProperties2} by GAP.
Here is a sample source code for GAP (Version 4.10.2)
which performs necessary calculation.

\lstset{
  language=GAP,
  breaklines = true,
  basicstyle=\ttfamily\scriptsize,
  breakindent = 10pt,
  numbers=left,
  numberstyle=\footnotesize
}

\begin{lstlisting}
G := PSp(6,2);;
C := List(ConjugacyClassesMaximalSubgroups(G), Representative);;
# G = Sp_6(F_2) has 8 maximal subgroups, with order
# 51840 40320 23040 12096 10752 4608 4320 1512

U28 := C[1];;
if Size(G)/Size(U28) = 28 then Print("OK. U28 has index 28.\n");
  else Print("Error: U28 does not have index 28.\n"); fi;

U36 := C[2];;
if Size(G)/Size(U36) = 36 then Print("OK. U36 has index 36.\n");
  else Print("Error: U36 does not have index 36.\n"); fi;

U63 := C[3];;
if Size(G)/Size(U63) = 63 then Print("OK. U63 has index 63.\n");
  else Print("Error: U63 does not have index 63.\n"); fi;

count_solvable_32 := function(GroupList)
  local count_solvable, count_32;
  count_solvable := 0;
  count_32 := 0;
  for K in GroupList do
    if IsSolvable(K) = true
      then count_solvable := count_solvable + 1; fi;
    if IsomorphismGroups(K, ElementaryAbelianGroup(32)) <> fail
      then count_32 := count_32 + 1; fi;
  od;
  Print("The number of solvable subgroups: ",count_solvable,"\n");
  Print("The number of subgroups isomorphic to (F_2)^5: ",count_32,"\n");
end;

count_cond := function(GroupList)
  local count_plus, count_plus_solvable, count_plus_32,
    count_minus, count_minus_solvable, count_minus_32,
    count_both, count_both_solvable, count_both_32;
  count_plus := 0;
  count_plus_solvable := 0;
  count_plus_32 := 0;
  count_minus := 0;
  count_minus_solvable := 0;
  count_minus_32 := 0;
  count_both := 0;
  count_both_solvable := 0;
  count_both_32 := 0;

  for K in GroupList do
    # Condition +
    A := Action(K, RightCosets(G, U36), OnRight);
    cond_plus := 0;
    if NrMovedPoints(A) < 36 then cond_plus := 1;
    else
      for x in A do if NrMovedPoints(x) = 36 then cond_plus := 1; fi; od;
    fi;

    # Condition -
    A := Action(K, RightCosets(G, U28), OnRight);
    cond_minus := 0;
    if NrMovedPoints(A) < 28 then cond_minus := 1;
    else
      for x in A do if NrMovedPoints(x) = 28 then cond_minus := 1; fi; od;
    fi;

    if cond_plus = 0 then
      count_plus := count_plus + 1;
      if IsSolvable(K) = true then
        count_plus_solvable := count_plus_solvable + 1; fi;
      if IsomorphismGroups(K, ElementaryAbelianGroup(32)) <> fail then
        count_plus_32 := count_plus_32 + 1; fi;
    fi;

    if cond_minus = 0 then
      count_minus := count_minus + 1;
      if IsSolvable(K) = true then
        count_minus_solvable := count_minus_solvable + 1; fi;
      if IsomorphismGroups(K, ElementaryAbelianGroup(32)) <> fail then
        count_minus_32 := count_minus_32 + 1; fi;
    fi;

    if cond_plus = 0 and cond_minus = 0 then
      count_both := count_both + 1;
      if IsSolvable(K) = true then
        count_both_solvable := count_both_solvable + 1; fi;
      if IsomorphismGroups(K, ElementaryAbelianGroup(32)) <> fail then
        count_both_32 := count_both_32 + 1; fi;
    fi;
  od;

  Print("Condition plus: ",count_plus,"\n");
  Print("Condition plus (solvable): ",count_plus_solvable,"\n");
  Print("Condition plus ((F_2)^5): ",count_plus_32,"\n");

  Print("Condition minus: ",count_minus,"\n");
  Print("Condition minus (solvable): ",count_minus_solvable,"\n");
  Print("Condition minus ((F_2)^5): ",count_minus_32,"\n");

  Print("Both: ",count_both,"\n");
  Print("Both (solvable): ",count_both_solvable,"\n");
  Print("Both ((F_2)^5): ",count_both_32,"\n\n");
end;

GroupList := List(ConjugacyClassesSubgroups(G), Representative);;
Print("Sp(6,2) has ",Size(GroupList)," subgroups, up to conjugation.\n");
count_solvable_32(GroupList);
count_cond(GroupList);

GroupList := List(ConjugacyClassesSubgroups(U28), Representative);;
Print("U28 has ",Size(GroupList)," subgroups, up to conjugation.\n");
count_solvable_32(GroupList);
count_cond(GroupList);

GroupList := List(ConjugacyClassesSubgroups(U36), Representative);;
Print("U36 has ",Size(GroupList)," subgroups, up to conjugation.\n");
count_solvable_32(GroupList);
count_cond(GroupList);

GroupList := List(ConjugacyClassesSubgroups(U63), Representative);;
Print("U63 has ",Size(GroupList)," subgroups, up to conjugation.\n");
count_solvable_32(GroupList);
count_cond(GroupList);
\end{lstlisting}

\vspace{0.1in}

If the above code is executed successfully,
it outputs as follows.

\begin{lstlisting}
OK. U28 has index 28.
OK. U36 has index 36.
OK. U63 has index 63.
Sp(6,2) has 1369 subgroups, up to conjugation.
The number of solvable subgroups: 1301
The number of subgroups isomorphic to (F_2)^5: 6
Condition plus: 411
Condition plus (solvable): 399
Condition plus ((F_2)^5): 6
Condition minus: 371
Condition minus (solvable): 359
Condition minus ((F_2)^5): 6
Both: 240
Both (solvable): 228
Both ((F_2)^5): 6
U28 has 350 subgroups, up to conjugation.
The number of solvable subgroups: 331
The number of subgroups isomorphic to (F_2)^5: 0
Condition plus: 22
Condition plus (solvable): 22
Condition plus ((F_2)^5): 0
Condition minus: 0
Condition minus (solvable): 0
Condition minus ((F_2)^5): 0
Both: 0
Both (solvable): 0
Both ((F_2)^5): 0
U36 has 296 subgroups, up to conjugation.
The number of solvable subgroups: 268
The number of subgroups isomorphic to (F_2)^5: 0
Condition plus: 0
Condition plus (solvable): 0
Condition plus ((F_2)^5): 0
Condition minus: 35
Condition minus (solvable): 35
Condition minus ((F_2)^5): 0
Both: 0
Both (solvable): 0
Both ((F_2)^5): 0
U63 has 1916 subgroups, up to conjugation.
The number of solvable subgroups: 1880
The number of subgroups isomorphic to (F_2)^5: 13
Condition plus: 856
Condition plus (solvable): 844
Condition plus ((F_2)^5): 13
Condition minus: 711
Condition minus (solvable): 699
Condition minus ((F_2)^5): 13
Both: 548
Both (solvable): 536
Both ((F_2)^5): 13
\end{lstlisting}

\vspace{0.1in}

Here is a sample source code for SageMath (Version 8.9)
which calculates quartics by the algorithm described in Section \ref{Sect:Algorithm}.
It also checks the smoothness of the output.
The example in Section \ref{Section:Example} is calculated by
this code.

\begin{lstlisting}
P.<X, Y, Z, S, T> = PolynomialRing(QQ)

b1, b2, b3, b4, b5 = -1, 17, 89, 257, 769
a1, a2, a3, a4, a5, u = b1*b5, b2*b4, b3, b4, b5, -b4^(-1)*b5^(-1)

a6 = a1*a2*a3*a4*a5*u^2
c = a1*a2*a3*a4*a5*u

F = (a1*S-T)*(a2*S-T)*(a3*S-T)*(a4*S-T)*(a5*S-T)*(a6*S-T)
g0 = (-F.coefficient({S:1, T:5}) -1)/2
g1 = (-F.coefficient({S:2, T:4}) + g0^2)/2
g2 = c
g = g0*T^2 + g1*S*T + g2*S^2
h = (T^5*S - (T^3 - g*S)^2 + F)/(S^3*T)
#print(factor(g))
#print(factor(h))

M = matrix([[-S*T+T^2, S*T, g], [S*T, S^2, T^2], [g, T^2, h]])
#print(det(M).factor())  ## check det(M) = -F(S,T)
v = matrix([X,Y,Z])
Biquad = (v*M*v.transpose())[0][0]
Biquad = P(Biquad)
#print(factor(Biquad))

q0 = Biquad.coefficient({S:2, T:0})
q1 = Biquad.coefficient({S:1, T:1})
q2 = Biquad.coefficient({S:0, T:2})
Quart = q1^2-4*q0*q2
#print(factor(Quart))

PP.<X, Y, Z> = ProjectiveSpace(QQ, 2)
R = PP.subscheme([Quart])
R.is_smooth()
\end{lstlisting}

\section{Complexity of our algorithm}
\label{Appendix:Complexity}

In Section \ref{Sect:Algorithm},
we give an algorithm to obtain smooth quartics over
global fields failing the local-global property for bitangents explicitly.
Our algorithm is fast, in practice.
But it seems difficult to analyze the time complexity by the following reasons.

For a given global field $K$, the hardest step in our algorithm is to find a Galois extension $L/K$
with $\Gal(L/K) \cong \F_2^{\oplus 5}$
such that every decomposition group is cyclic; see Step 1 in Section \ref{Sect:Algorithm}.
(For example, when $K = \Q$, the extension
\[ \Q\big(\,\sqrt{-1},\,\sqrt{17},\,\sqrt{89},\,\sqrt{257},\,\sqrt{769} \,\big)/\Q \]
satisfies this condition.
From this, we get the smooth quartic over $\Q$ in Section \ref{Section:Example}.)

In practice, it is easy to find an extension $L/K$.
But it is not easy to give a sharp estimate of the time complexity.
Assume the generalized Riemann hypothesis (GRH) for simplicity.
Then, for a finite Galois extension $M/K$ and a conjugacy class $C \subset \Gal(M/K)$,
by the effective Chebotarev density theorem,
we can find a prime $\mathfrak{p}$ of $K$
such that the Frobenius element at $\mathfrak{p}$
belongs to $C$ and
\[ N_{K/\Q}(\mathfrak{p}) \leq c (\log d_M)^2 (\log \log d_M)^4. \]
Here $d_M$ is the absolute value of the discriminant of $M$,
and $c$ is an effectively computable absolute constant;
see \cite[Corollary 1.2]{LagariasOdlyzko}.
From this, we can estimate
the time complexity to find an extension $L/K$
by a polynomial function with respect to $\log d_K$.
(Without assuming GRH, we can estimate
the time complexity by an exponential function with respect to $\log d_K$.)

Once the extension $L/K$ is found,
the rest is a simple calculation of the determinant,
which can be done in $O(1)$ arithmetic operations in the number field $K$.
(See the sample source code for SageMath in Appendix \ref{Appendix:SourceCode}.)
However, note that the quartic obtained by our algorithm
might be singular.
If it is singular, we may take other parameters randomly and calculate the quartic again.
The quartic obtained by our algorithm is smooth if parameters are chosen generically; see \cite[Proposition 3.5]{ElsenhansJahnel:SteinerHexad}.

\subsection*{Acknowledgements}

The authors would like to thank J\"org Jahnel and
Daniel Loughran for explaining how to construct quartics failing
the local-global property for bitangents by the results in \cite{JahnelLoughran};
see Remark \ref{Remark:JahnelLoughranl}.
The authors would like to thank the referee for comments and
constructive advice.
Y.\ I.'s work was supported by JSPS KAKENHI Grant Number 16K17572.
T.\ I.'s work was supported by JSPS KAKENHI Grant Number 20674001 and 26800013.
T.\ O.'s work was supported by JSPS KAKENHI Grant Number 26800011 and 18H05233.
T.\ T.'s work was supported by JSPS KAKENHI Grant Number 17H02835.
Y.\ U.'s work was supported by JSPS KAKENHI Grant Number 20K03517.
This work was supported by the Sumitomo Foundation FY2018 Grant
for Basic Science Research Projects (Grant Number 180044).
Most of calculations were done
with the aid of GAP \cite{GAP}, Maxima \cite{Maxima}, Sage \cite{Sage}.

\end{document}